\documentclass[11pt]{article}
\usepackage{amsmath,amssymb,amsfonts,amsthm,bm}
\usepackage{mathrsfs}
\usepackage{indentfirst}
\usepackage[pdftex]{color,graphicx}
\usepackage{subfig}
\usepackage[pdftex,bookmarks,unicode,colorlinks]{hyperref}
\usepackage{enumerate} 
\usepackage[numbers]{natbib} 
\usepackage{yfonts} 
\usepackage{caption} 
\usepackage{mathtools}
\usepackage{chngcntr}
\usepackage{algorithmicx,algorithm} 
\usepackage{lineno}

\theoremstyle{plain}

\theoremstyle{remark}

\theoremstyle{definition}

\numberwithin{equation}{section} 

\begin{document}
\title{\bf\Large The Most Likely Transition Path  \\ for a Class of Distribution-Dependent Stochastic Systems }
\author{\bf\normalsize{
Wei Wei$^{1,}$\footnotemark[2]
~and Jianyu Hu$^{2,}$\footnotemark[1]
}\\[10pt]
\footnotesize{$^1$Institute of Natural Sciences,
Shanghai Jiao Tong University,} \\
\footnotesize{Shanghai 200240, China} \\
\footnotesize{$^2$Division of Mathematical Sciences,
Nanyang Technological University,} \\
\footnotesize{21 Nanyang Link, 637371, Singapore.} 
}

\footnotetext[2]{Email: \texttt{weiw\_sjtu@sjtu.edu.cn}}
\footnotetext[1]{Email: \texttt{jianyu.hu@ntu.edu.sg}}
\footnotetext[1]{is the corresponding author}

\date{}
\maketitle
\vspace{-0.3in}

\begin{abstract}
    Distribution-dependent stochastic dynamical systems arise widely in engineering and science. We consider a class of such systems which model the limit behaviors of interacting particles moving in a vector field with random fluctuations. We aim to examine the most likely transition path between equilibrium stable states of the vector field. In the small noise regime, the action functional does not involve the solution of the skeleton equation which describes the unperturbed deterministic flow of the vector field shifted by the interaction at zero distance. As a result, we are led to study the most likely transition path for a stochastic differential equation without distribution dependency. This enables the computation of the most likely transition path for these distribution-dependent stochastic dynamical systems by the adaptive minimum action method and we illustrate our approach in two examples.

    \bigskip
 
\end{abstract}

\section{Introduction}
Large systems of interacting particles are now quite common in many fields. Particles can represent ions and electrons in plasmas in the context of physics \cite{Pope2000}. In the bio-science, they can model the motion of cells or bacteria and in economics or social sciences, they represent ``agents" \cite{Carmona17a, Carmona17b}. But they are analytically complex and computationally expensive as the number of particles is very large. Instead, a corresponding mean-field limit equation is often an alternative model to investigate. There are rich literatures on deriving such mean-field limit equations; see, for example, \cite{Bolley07, CS07, Jabin2014}.  Interactions in such equations often appear as the distribution in the drift term, which makes them known as distribution-dependent equations. The McKean-Vlasov stochastic differential equation was first studied in \cite{MKe} and is a widely applicable one among these distribution-dependent equations. 

The motion of the particles modeled by the McKean-Vlasov stochastic differential equation can be simply described as follows: particles interacting with each other wander under random fluctuations in a landscape, which is formed by a vector field. We are interested in how such a particle transfers from one equilibrium stable state of the landscape to another equilibrium stable state. Large deviations theory is a useful tool to deal with this kind of problem. Roughly speaking, the large deviations theory characterizes how small the probability of a rare event is via an action functional. It also extends the concept of potential to non-gradient systems as quasi-potential, which is widely used in investigating the long time behaviors of a system under small random perturbations. More details of the large deviations theory may be found in \cite{DZ98,FW12}. 

There are many large deviations results investigating the McKean–Vlasov stochastic differential equations.  Herrmann et al.~\cite{Imkeller2008} derived a large deviation result for a one-dimensional McKean-Vlasov stochastic equation. Based on this, they investigated the influence of interactions on exit times and exit locations. Recently, Reis et al.~\cite{Reis2018} obtained large deviations results for a more general type of $n$-dimensional McKean-Vlasov stochastic equations with Gaussian noise. Liu et al.\cite{liu2020} extended the large deviations results to the McKean–Vlasov stochastic equations with jumps. Although these equations contain the distribution of the solution in the drift term, the action functionals they obtained are simplified to only including a Dirac measure. Therefore, numerical methods built on the large deviations theory can be applied to investigate the dynamical behaviors of such McKean-Vlasov stochastic equations.

The minimum action method(MAM) is a widely used numerical method in finding the most likely transition path. It was proposed to compute the most likely transition path for stochastic differential equations with multiple metastable states~\cite{weinan2004}. The basic idea of MAM is to find the minimizer of the associated action functional and take this minimizer as the desired most likely transition path. MAM was improved by a reparametrization \cite{heymann2008} to a fixed finite time scale and numerically solving an Euler-Lagrange equation to find such a minimizer. A moving mesh method was used in \cite{sun2018, zhou2008} to enhance the performance of MAM. Authors in \cite{wan2018} employed finite elements approximation to find the minimizer of the associated action functional.

In this paper, we calculate the action functionals for a class of McKean–Vlasov stochastic differential equations. We find that the action functional, which characterizes the probability of rare events, usually involves a `modified deterministic flow', which is the flow of the vector field shifted by the interaction at zero distance. When we examine the most likely transition path between equilibrium stable states of the vector field, we find that the action functional does not involve this modified  deterministic flow. A numerical method based on the adaptive minimum action method \cite{sun2018,zhou2008} is used to calculate the most likely transition path. 

The paper is arranged as follows. In section \ref{bkg}, we briefly introduce the background of the McKean–Vlasov stochastic equation and the associated large deviations theory. Based on the large deviations theory, we analyze the transition behaviors of particles modeled by the McKean–Vlasov stochastic equation. In section \ref{numeri}, we will introduce the numerical method for computing the most likely transition path and present the most likely transition paths for a bi-stable system under two different interactions.


\section{Backgrounds}\label{bkg}

We consider the following $d$-dimensional McKean-Vlasov equation with $X^\epsilon_0=x_1$ , 
\begin{equation}\label{dde}
    d X_{t}^{\epsilon}=V\left(X_{t}^{\epsilon}\right) d t- F*u_t^\epsilon(X_t^\epsilon)d t+\sqrt{\epsilon}  d B_t,
\end{equation}
where $u^\epsilon_t$ stands for the law of $X^\epsilon_t$ and solution $X^\epsilon$ is a stochastic process on the given probability space $(\Omega, \mathcal{F}, \mathbb{P})$. The symbol ``$*$'' is the usual convolution operator that, for $\omega_0 \in \Omega$
\begin{align*}
    F *u_s^\epsilon(X_s^\epsilon(\omega_0))& = \int_{\mathbb{R}^d}F(X_s^\epsilon(\omega_0)-x) d u_s^\epsilon(dx) \\
    &=\int_{\omega \in \Omega} F \left( X_s^\epsilon(\omega_0)-X_s^\epsilon (\omega)\right) d \mathbb{P}(\omega) = \mathbb{E}\left(F( X_s^\epsilon(\omega_0)-X_s^\epsilon )\right).
\end{align*} 
Equation \eqref{dde} can be regarded as the mean-field limit of the following system of many particles as $N \to \infty$, which is known as the \emph{propagation of chaos}.

Let $X^{\varepsilon}_{i, N}(0)=x_1$ and 
\begin{equation}\label{mfe}
    d X^{\varepsilon}_{i, N}(t)=V\left(X^{\varepsilon}_{i, N}(t)\right) d t-\frac{1}{N} \sum_{j=1}^{N} F\left(X^{\varepsilon}_{i, N}(t)-X^{\varepsilon}_{j, N}(t)\right) dt+\sqrt{\varepsilon} d B_{i, N}(t),
\end{equation}
where the Brownian motions $B_{i,N}(t)$ are independent and $F$ represents the interaction between particles. 

The solution $X^{\varepsilon}_{i, N}$ is regarded as the motion of the $N$-th particle in a random environment that interacts with all other particles therein. The solution $X^\epsilon$ can be seen as the motion of a particle in a macroscopic system as $N$ in equation (\ref{mfe}) is the same magnitude as Avogadro's constant. The motion of $X^\epsilon$ is generated by the vector field $V$, the interaction $F$, and the random fluctuation $B_t$. 

We are interested in the transition behaviors of $X^\epsilon$ when the vector field $V$ admits multiple equilibrium states, that is $V=0$ at these states. Will these states become the metastable states of $X^\epsilon$ or will the interaction $F$ create a different metastable state? What does the most likely transition path look like? The large deviations theory has shown its power in tackling such kinds of problems. We will briefly introduce the large deviations theory in the following subsection.

\subsection{The large deviations theory}
 
Let $C([0,T])$ denote the space of continuous functions on the interval $[0,T]$ taking values in $\mathbb{R}^d$, equipped with the uniform norm. Roughly speaking, if the solution $X^\epsilon$ to (\ref{dde}) satisfies the large deviations principle with an action functional $I_T$ on $C([0,T])$, we can estimate the probability that $X^\epsilon$ stays in a neighborhood of a path $\varphi$ as follows, 
\begin{equation*}
    \mathbb{P}\{ \| X^\epsilon-\varphi\|<\delta\} \sim \exp \left(-\frac{1}{\epsilon}I_T(\varphi)\right), \quad \textrm{as } \epsilon \to 0.
\end{equation*}
The above asymptotic results can be used to estimate the probability of events associated with $X^\epsilon$ as $\epsilon$ tends to $0$. To find the most likely transition path of $X^\epsilon$ from point $x_1$ to point $x_2$ within time $T$, we only need to solve the following constrained minimization problem. That is 
\begin{equation}\label{mini}
    I_T(\varphi^*) =  \min_{\varphi} I_T(\varphi),
\end{equation}
where the minimization is constrained by 
\begin{equation*}
    \varphi(0) = x_1, \quad \varphi(T) = x_2, \quad \varphi \in C([0,T]).
\end{equation*}
The function $\varphi^*_T$ gives the most likely transition path within time $T$. The infimum value of equation \eqref{mini} over $T\in [0,\infty]$ is called quasi-potential $\tilde{V}(x_1,x_2)$ between $x_1$ and $x_2$, that is 
\begin{equation}\label{MTP}
    \tilde{V}(x_1,x_2) = \inf_{T>0}\inf_{\varphi \in \bar{C}_{x_1}^{x_2}([0,T])} I_T(\varphi),
\end{equation}
where $\bar{C}_{x_1}^{x_2}(0,T)$ is the space of all absolutely continuous functions that start at $x_2$ and end at $x_2$. Suppose that there exists a path $\tilde{\varphi}$ and time $\tilde{T} \in [0,\infty]$ satisfying,
\begin{equation*}
    I_{\tilde{T}}(\tilde{\varphi}) = \tilde{V}(x_1,x_2).
\end{equation*}
The path $\tilde{\varphi}$ is the most likely transition path from point $x_1$ to $x_2$ that we are looking for.


As is shown above, the action functional $I_T$ governs the dynamical behavior of a stochastic system in the context of the large deviations theory. With results derived in \cite{Reis2018}, we are able to calculate the action functional $I_T$ for our McKean–Vlasov equation (\ref{dde}). 

Suppose the coefficient $b(x,t,\mu_t):=V(x)-F*\mu_t (x) $ satisfies \cite[Assumption 3.2]{Reis2018}, where $\mu_t$ is the law of a random variable. Let $\eta$ be the solution to the following ordinary differential equation,
\begin{equation}\label{eta}
    \dot{\eta}(t)=V(\eta(t))-F * \delta_{\eta(t)}(\eta(t))=V(\eta(t))-F(0), \quad \eta(0)=x_1.
\end{equation} 
Equation (\ref{eta}) is often regarded as the skeleton equation. According to \cite[Theorem 4.4]{Reis2018}, we know that $X^\epsilon$ satisfies the large deviations principle in $C([0,T])$ equipped with the topology of the uniform norm with the action functional 
\begin{equation}\label{Act}
    I^{x_1}_T(\varphi)=\frac{1}{2} \int_0^T \left|\dot{\varphi}(s)-V(\varphi(s))+F(\varphi(s)-\eta(s))\right|^2 ds,
\end{equation}
where $\varphi$ is an absolutely continuous function.
$X_t^\epsilon$ shares the same action functional with the following corresponding stochastic differential equation if it satisfies the large deviations principle,
\begin{equation}\label{dsde}
    dZ^\epsilon_t=\left(V(Z^\epsilon_t)-F(Z^\epsilon_t-\eta(t))\right)dt+\sqrt{\epsilon}dB_t.
\end{equation}

Note that the drift term of \eqref{dsde} may involve time $t$ explicitly.

\subsection{Dynamical behaviors of the McKean–Vlasov stochastic differential equation}
We are interested in the transition behavior of the McKean–Vlasov equation between the equilibrium stable states of the vector field $V$. Based on the large deviation theory, we will analyze the dynamical behavior of a class of McKean–Vlasov stochastic differential equations in this subsection.

Let $x_1$ and $x_2$ be two equilibrium stable states of the vector field $V$, that is $V(x_1)=V(x_2)=\mathbf{0}$. If the interaction $F$ is $\mathbf{0}$ at point $\mathbf{0}$, equation \eqref{eta} starting at $x_1$ reduces to 
\begin{equation}\label{exeta}
    \dot{\eta}(t)=V(\eta(t)), \quad \eta(0)=x_1.
\end{equation}
 When we look for the most likely transition path from $x_1$ to $x_2$, the solution $\eta$ stays at $x_1$. Then, the action functional in this case is 
\begin{equation*}
    I^{x_1}_T(\varphi)=\frac{1}{2} \int_0^T \left|\dot{\varphi}(s)-V(\varphi(s))+F(\varphi(s)-x_1)\right|^2 ds  = \int_0^T |\dot{\varphi}-\bar{b}(\varphi(s))|^2ds.
\end{equation*}
The corresponding stochastic differential equation (\ref{dsde}) then becomes 
\begin{equation}\label{Csde}
    dZ^\epsilon_t=\left(V(Z^\epsilon_t)-F(Z^\epsilon_t-x_1)\right)dt+\sqrt{\epsilon}dB_t=\bar{b}(Z_t^\epsilon)dt+\sqrt{\epsilon}dB_t, \quad Z_0^\epsilon = x_1.
\end{equation}
For the solution $\xi$ to the following equation,
\begin{equation}\label{CsdeS}
    \dot{\xi}(t)=V(\xi(t))-F(\xi(t)-x_1), \quad  0\leq t < \infty,
\end{equation}
we know that $I_t^{x_1}(\xi)=0$. Because we can choose $\xi(0)$ as the starting point and $\xi(t)$ as the end point. This is also the minimum value of the action functional $I_t^{x_1}$ on absolutely continuous functions from $\xi(0)$ to $\xi(t)$. So the quasi-potential between every two different points on the path $\xi$ is $0$. Then the trajectories of $\xi$ starting at different points form a set of equipotential lines.
In addition, we assume that the equilibrium stable states $x_1$ and $x_2$ are also asymptotic stable states of the system (\ref{exeta}). Then for the following system
\begin{equation} \label{SSde}
    dY_t^\epsilon=V(Y_t^\epsilon)dt+\sqrt{\epsilon} dB_t,
\end{equation}
$x_1$ and $x_2$ become two metastable states. It means that $Y^\epsilon$ starting near one of these two states, will stay a fairly long time near this state. But for the solution $Z^\epsilon$ to the corresponding stochastic differential equation (\ref{Csde}), $x_1$ and $x_2$ may not be metastable states due to the presence of the interaction term $F$. Instead, we need to find the equilibrium states of equation \eqref{CsdeS} to determine the metastable states of $Z^\epsilon$. As is shown in the next section, the most likely transition path will pass by a saddle point and be attracted to a ``interaction-generated" metastable state of $Z^\epsilon$.

Note that the integrand $\bar{b}$ does not involve time $t$ explicitly. We can use a reparametrization method to find the most likely transition path over time $T \in [0,\infty]$, such as the geometric minimum action method \cite{heymann2008}. For the cases that the drift term $b$ in equation (\ref{dsde}) involves time $t$ explicitly, the adaptive minimum action method \cite{sun2018} can be used to compute the most likely transition path.


\section{Compute the most likely transition path for the McKean-Vlasov stochastic differential equation}\label{numeri}

In this section, we will introduce the numerical method based on the adaptive minimum action method \cite{sun2018} for computing the most likely transition path for the distribution-dependent system \eqref{dde}.

The idea of the adaptive minimum action method(aMAM) is to find the minimizer of \eqref{MTP}, by numerically solving an associated Euler-Lagrange equation with a large enough time $T$ for constrained minimization problem \eqref{mini}. The moving mesh strategy is used to choose a proper mesh grid via a monitor function in every iteration. The associated Euler-Lagrange equation is then solved on such mesh grids. Choose a large enough $T$ and the numerical solution is the most likely transition path we are looking for.

To solve the constrained minimization problem \eqref{mini}, we will numerically solve the associated Euler-Lagrange equation whose solution is regarded as the steady state of the gradient flow $\delta I_T / \delta \varphi$ Specifically, the Euler-Lagrange equation for action functional \eqref{Act} is 
\begin{equation}\label{EL}
    \left\{ \begin{array}{l} 
         \varphi_{tt}-\nabla_x b(\varphi,\eta) \varphi_t+\big(\nabla_x b(\varphi,\eta)\big)^\top (\varphi_t-b(\varphi,\eta))-\nabla_y b(\varphi,\eta) b(\eta,\eta)=0, \\
         \varphi(0)=x_1, \quad \varphi(T)=x_2
    \end{array}\right.
\end{equation}
where $\eta $ is the solution of equation \eqref{eta} and $b(x,y)=V(x)-F(x-y)$. 

The adaptive mesh grids for solving equation \eqref{EL} are determined by a monitor function. In this paper, we adopt the monitor function in \cite{sun2018} and it is given by $w(s)=|b(\varphi(s),\eta(s))|^r/C$ with $C=\int_0^T|b(\varphi(s),\eta(s))|^r d s$. The new variable $\alpha$ is then given by $\alpha(s)=\int_0^s w(\tau) d\tau $. Denote $\tilde{\varphi}(\alpha(s))=\varphi(s)$. Then, with respect to parameter $\alpha$, the Euler-Lagrange equation \eqref{EL} becomes 
\begin{equation*}
    \left\{ \begin{array}{l} 
        w^2\tilde{\varphi}^{\prime \prime}+\biggl(\big(\nabla_x b(\tilde{\varphi},\eta)\big)^\top-\nabla_x b(\tilde{\varphi},\eta)+w_\alpha \biggr)w\tilde{\varphi}^\prime+ \nabla_x b(\tilde{\varphi},\eta) b(\tilde{\varphi},\eta)-\big(\nabla_y b(\tilde{\varphi},\eta)\big)^\top b(\eta,\eta) =0,\\
        \tilde{\varphi}(0)=x_1, \quad \tilde{\varphi}(T)=x_2,
    \end{array}\right.
\end{equation*}
where $\tilde{\varphi}^\prime$ and $w_\alpha$ is the derivative with respect to parameter $\alpha$.  The action functional \eqref{Act} then becomes 
\begin{equation*}
    I_T(\varphi)=I(\tilde{\varphi})=\frac{1}{2}\int_0^1 |\dot{\tilde{\varphi}}(\alpha)w(s)-V(\tilde{\varphi}(\alpha))+F(\tilde{\varphi}(\alpha)-\tilde{\eta}(\alpha))|^2 \frac{d \alpha}{w(s)}.
\end{equation*}
The algorithm is given in Algorithm \ref{alg}.

\begin{algorithm}[htbp]
    \caption{An algorithm for distribution-dependent equations}\label{alg}
        \begin{algorithmic}[1]
            \State At the k-th iteration, given $(t_i^k,\varphi_i^k,\eta_i^k)$ for $i=0,1,\cdots,N$, calculate $w_i^k=|b(\varphi_i^k,\eta_i^k)|^r/C_k$ with $C_k=\sum_i (t_{i+1}^k-t_i^k) |b\big((\varphi_{i+1}^k+\varphi_i^k)/2,(\eta_{i+1}^k+\eta_i^k)/2\big)|^r$.
            \State Calculate 
            \begin{equation*}
                \Delta \alpha_i^k=\frac{(t_{i+1}^k-t_i^k) \left|b\left((\varphi_{i+1}^k+\varphi_i^k)/2,(\eta_{i+1}^k+\eta_i^k)/2\right) \right|^r}{C_k}.
            \end{equation*}
            \State Set $\Delta \alpha=1/N$ and $\alpha_{i+1}-\alpha_i=\Delta \alpha$ for $i=0,\cdots, N $ and $\alpha_0=0$.
            \State Interpolate $(\alpha_i^k,\varphi_i^k)$ to get $(\alpha_i,\bar{\varphi}_i^k)$, and interpolate $(\alpha_i^k, t_i^k)$ to get $(\alpha_i,t_i^{k+1},\eta_i^{k+1})$, where $\eta_i^k$ is generated by the Euler scheme at time $t_i^k$.
            \State Let $\{\bar{\varphi}_i\}_{i=1}^{N+1}$ be the solution of the following linear equations 
            \begin{equation*}
                \left\{ 
                    \begin{aligned}
                        \frac{\bar{\varphi}_i-\bar{\varphi}_i^k}{\Delta \tau}  = & (w_i^k)^2 \frac{\bar{\varphi}_{i+1}-2\bar{\varphi}_i+\bar{\varphi}_{i-1}}{\Delta \alpha^2} + \bigg(\nabla_y b(\bar{\varphi}_i^k,\eta_i^{k+1})^\top - \nabla_y b(\bar{\varphi}_i^k,\eta_i^{k+1})+(w_\alpha)_i^k \textrm{Id}\bigg) w_i^k \bar{\varphi}_i^{\prime k}\\
                        & -\nabla_y b(\bar{\varphi}_i^k,\eta_i^{k+1}) b(\bar{\varphi}_i^k,\eta_i^{k+1})-\nabla_x b(\bar{\varphi}_i^k,\eta_i^{k+1})^\top b(\eta_i^{k+1},\eta_i^{k+1}) \\
                        \bar{\varphi}_0 = x_1, \quad & \bar{\varphi}_N=x_2.
                    \end{aligned}\right.
            \end{equation*}
            where $\bar{\varphi}_i^{\prime k}=(\bar{\varphi}_{i+1}^k-\bar{\varphi}_{i-1}^k)/(2 \Delta \alpha)$ and $(w_\alpha)_i^k=(w_{i+1}^k-w_{i-1}^k)/(2\Delta \alpha)$
            \State Repeat step 1-5 until a stopping criterion is fulfilled.
        \end{algorithmic}
\end{algorithm}

We will compute the most likely transition path for the Maier-Stein model with different interaction terms. That is 
\begin{equation*}
    d X_{t}^{\epsilon}=V\left(X_{t}^{\epsilon}\right) d t- F*u_t^\epsilon(X_t^\epsilon)d t+\sqrt{\epsilon}d B_t, \quad X_0^\epsilon=(-1,0),
\end{equation*}
where
\[
    V(u,v)=\left(\begin{array}{c}
    u-u^{3}-\beta u v^{2} \\-\left(1+u^{2}\right) v\end{array}
    \right), \quad \beta=10.
\]


\begin{figure}[htbp]
    \centering  
    \includegraphics[width=0.7\textwidth]{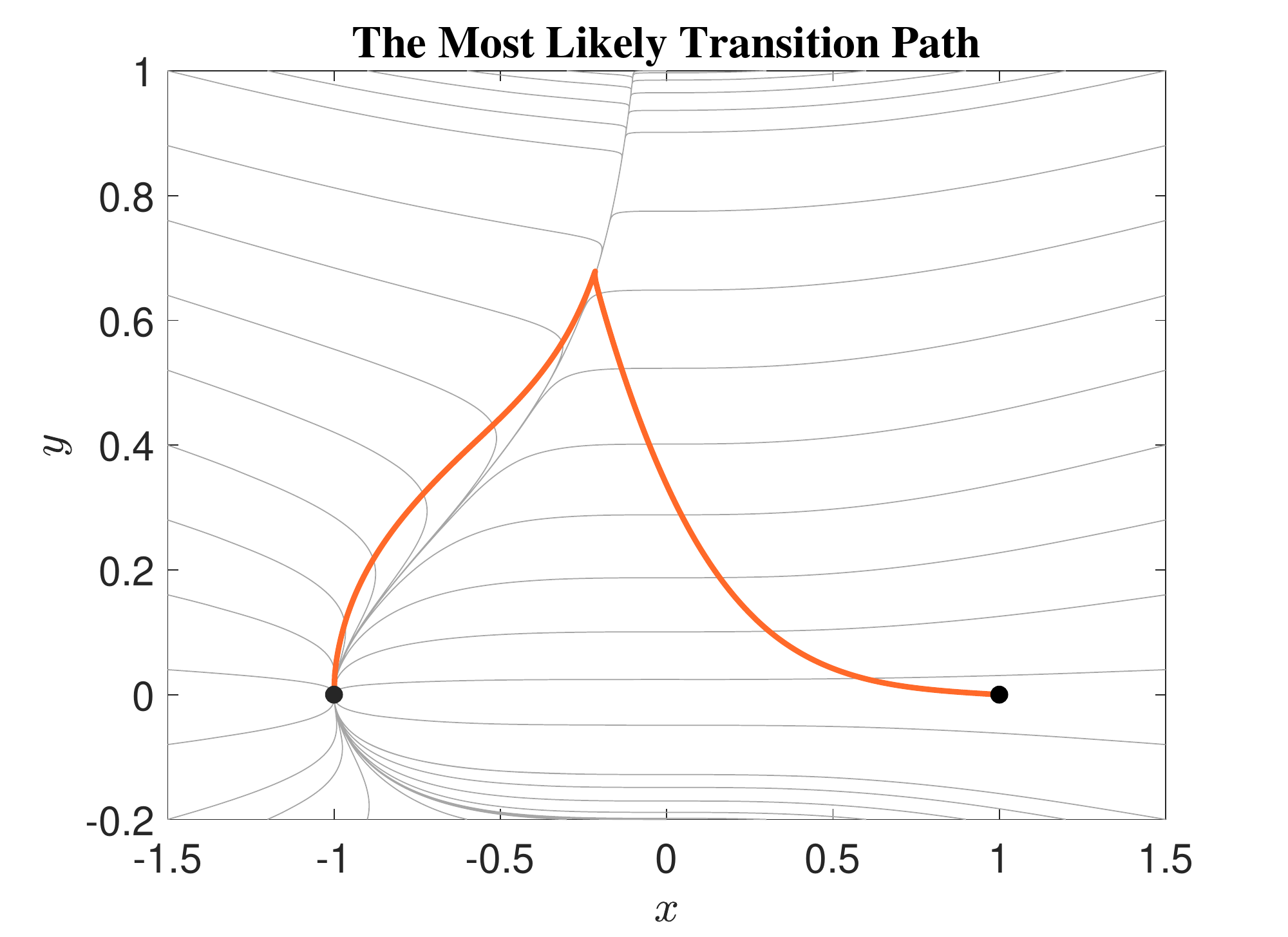}
    \caption{The most likely transition path for the Maier-Stein model with interaction $F_L$. The orange line is the most likely transition path connecting points $(-1,0)$ and $(1,0)$. Grey lines are the equipotential lines of the associated quasi-potential. We set the number of points $N = 200$, iteration $K =5000$ and time $T =20$. The initial path is chosen as $y=-0.5x^2+0.5$.}\label{Line}
\end{figure}

In Figure~\ref{Line}, we calculate the most likely transition path connecting two equilibrium stable states $(-1,0)$ and $(1,0)$, under the interaction $F=F_{L}$, which is 
\begin{equation*}
    F_{L}(u,v)=\left(\begin{array}{c}
        u \\ -v
    \end{array}\right).
\end{equation*}

The corresponding stochastic differential equation (\ref{Csde}) then becomes
\begin{equation*}
    dZ^\epsilon_t=b_L(Z_t^\epsilon)dt+\sqrt{\epsilon}dB_t, \quad Z_0^\epsilon = (-1,0),
\end{equation*}
where 
\begin{equation*}
    b_L(u,v)=\left( 
        \begin{array}{c}
            -u^3-\beta u v^2+1 \\
            -u^2 v
        \end{array}\right).
\end{equation*}
The action functional is 
\begin{equation*}
    I_L(\varphi)=\int_0^T |\dot{\varphi}-b_L(\varphi(s))|^2ds.
\end{equation*}
There is only one asymptotic stable point $(-1,0)$ for the following equation, 
\begin{equation*}
    \dot{\xi_L}(t)=b_L(\xi_L(t)).
\end{equation*}
It is also the metastable point of the corresponding stochastic differential equation. The grey lines in Figure~\ref{Line} are the orbits of $\xi_L$ starting at different points and these orbits form the equipotential lines of the quasi-potential. The interaction $F$ eliminates the metastable point $(1,0)$. The most likely transition path in this case is actually the most likely exit path from $(-1,0)$ to $(1,0)$. 

\begin{figure}[htbp]
    \centering
    \includegraphics[width=0.85\textwidth]{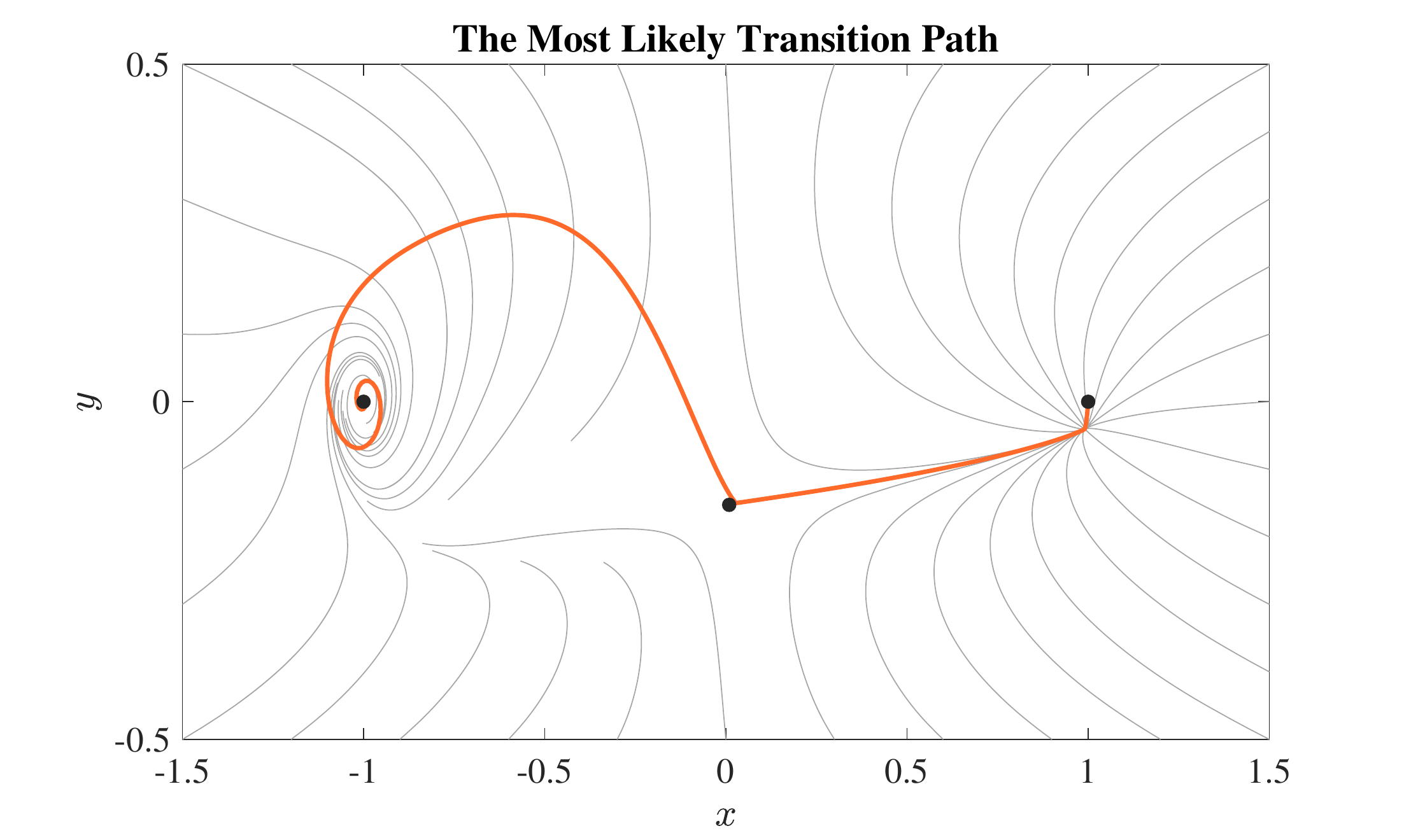}
    \caption{The most likely transition path for the Maier-Stein model with interaction $F_{BS}$. The orange line is the most likely transition path connecting points $(-1,0)$ and $(1,0)$. Grey lines are the equipotential lines of the associated quasi-potential.  The dark grey point is the saddle point of equation \eqref{eBS}. We set the number of points $N = 400$, iteration $K =15000$ and time $T =40$. The initial path is chosen as $y=-0.5x^2+0.5$. }\label{BS}
\end{figure}
In Figure~\ref{BS}, we calculate the most likely transition path under the interaction described by a modified Biot-Savart kernel which is not singular at the origin. Precisely, we take $F=F_{BS}$ to be
\begin{equation*}
    F_{BS}(u,v)=\frac{1}{2 \pi}\left( \begin{array}{c}
        \frac{-v}{u^2+v^2+\delta} \\
        \frac{u}{u^2+v^2+\delta}\end{array}  \right), \quad \delta=0.01.
\end{equation*}
The corresponding stochastic differential equation (\ref{Csde}) then becomes
\begin{equation*}
    dZ^\epsilon_t=b_{BS} (Z_t^\epsilon)dt+\sqrt{\epsilon}dB_t,\quad Z_0^\epsilon = (-1,0),
\end{equation*}
where 
\begin{equation*}
    b_{BS}(u,v)=\left( 
        \begin{array}{c}
            u-u^3-\beta u v^2+\frac{v}{2\pi((u+1)^2+v^2+\delta)} \\
            -(1+u^2) v-\frac{u+1}{2\pi((u+1)^2+v^2+\delta)}
        \end{array}\right).
\end{equation*}
The action functional is 
\begin{equation*}
    I_{BS}(\varphi)=\int_0^T |\dot{\varphi}-b_{BS}(\varphi(s))|^2ds.
\end{equation*}
For the presence of the interaction term $F_{BS}$, there are two asymptotic stable points $x_1=(-1,0)$, $x_2 \approx (1.14801,0.0318893)$ and a saddle point $x^* \approx (0.009145,-0.152703)$ for the following equation, 
\begin{equation}\label{eBS}
    \dot{\xi}_{BS}(t)=b_{BS}(\xi_{BS}(t)).
\end{equation}
The metastable points for the corresponding stochastic differential equation are then $x_1$ and $x_2$. The grey lines in Figure~\ref{BS} are the orbits of $\xi_{BS}$ starting at different points and these form the equipotential lines of the quasi-potential. The saddle point $x^*$ is the dark grey point in Figure~\ref{BS}. The interaction $F$ transports the saddle point $(0,0)$ and stable point $(1,0)$ to $x^*$ and $x_2$ respectively. Then, the most likely transition path will pass by the saddle point $x^*$ and be attracted to the neighborhood of $x_2$ as is shown in Figure~\ref{BS}.

\section{Conclusion}

We have investigated the most likely transition path for a class of distribution-dependent stochastic dynamical systems. In the small noise regime, this problem turns out to be equivalent to computing this transition path for a stochastic dynamical system {\it without} distribution-dependency. Therefore, we can apply the adaptive minimum action method to calculate the most likely transition paths, as demonstrated in several examples. 

\section*{Acknowledgements}
This work was partly supported by the NSFC grants 11771449.

\section*{Data Availability}
The data that support the findings of this study are openly available in GitHub: \url{https://github.com/JayWeiess/MVSDE-aMAM}


\begin{thebibliography}{10}

    \bibitem{Reis2018}
    G.~Reis, W.~Salkeld and J.~Tugaut.
    \newblock{Freidlin–Wentzell LDP in path space for McKean–Vlasov equations and the functional iterated logarithm law}.
    \newblock {\em The Annals of Applied Probability}, 29:1487 -- 1540, 2019

    \bibitem{DZ98}
    A.~Dembod and O.~Zeltouni.
    \newblock {\em Large Deviations Techniques and Applications}.
    \newblock Springer, Berlin, Heidelberg, 2010.
    
    \bibitem{FW12}
    M.~I.~Freidlin and A.~D.~Wentzell.
    \newblock {\em Random Perturbations of Dynamical Systems}.
    \newblock Springer, Berlin, Heidelberg, 2012

    \bibitem{Imkeller2008}
    S.~Herrmann, P.~Imkeller and D.~Peithmann.
    \newblock{Large deviations and a Kramers’ type law for self-stabilizing diffusions}.
    \newblock{\em The Annals of Applied Probability}, 18:1379 -- 1423, 2008.
    
    \bibitem{Bolley07}
    F.Bolley, A.~Guillin and C.~Villani.
    \newblock{Quantitative concentration inequalities for empirical measures on non-compact spaces}.
    \newblock{\em Probability Theory and Related Fields}, 137:541 -- 593, 2007.
    
    \bibitem{CS07}
    F.~Cucker and S.~Smale.
    \newblock{Emergent behavior in flocks}.
    \newblock{\em IEEE Transactions on Automatic Control}, 52:852-- 862, 2007.
    
    \bibitem{Jabin2014}
    P.E.~Jabin.
    \newblock{A review for the mean field limit for Vlasov equations}.
    \newblock{\em Kinetic and Related Models}, 7:661-- 711, 2014.

    \bibitem{MKe}
    H.~P.~McKean.
    \newblock {A class of Markov processes associated with nonlinear parabolic equations}.
    \newblock {\em Proceedings of the National Academy of Sciences}, 56: 1907 -- 1911, 1966.

    \bibitem{Carmona17a}
    R.~Carmonaand F.~Delarue.
    \newblock{ \em Probabilistic Theory of Mean Field Games with Applications I}. 
    \newblock Springer, Cham, 2017.
    

    \bibitem{Carmona17b}
    R.~Carmonaand F.~Delarue.
    \newblock{\em Probabilistic Theory of Mean Field Games with Applications II}.   
    \newblock Springer, Cham, 2017.

    \bibitem{Pope2000}
    S.~B.~Pope.
    \newblock{\em Turbulent Flows}.
    \newblock Cambridge University Press, 2000.
    
    \bibitem{liu2020}
    W.~Liu, Y.~Song, J.~Zhai and T.~Zhang.
    \newblock{Large and moderate deviation principles for McKean-Vlasov SDEs with jumps}.
    \newblock{\em arXiv preprint arXiv:2011.08403} 2020

    \bibitem{heymann2008}
    M.~Heymann and E.~Vanden-Eijnden.
    \newblock {The geometric minimum action method: A least action principle on the space of curves}.
    \newblock {\em Communications on Pure and Applied Mathematics: A Journal Issued by the Courant Institute of Mathematical Sciences}, 61: 1052--1117, 2008.

    \bibitem{sun2018}
    Y.~Sun and X.~Zhou.
    \newblock {An improved adaptive minimum action method for the calculation of transition path in non-gradient systems}.
    \newblock {\em Communications in Computational Physics}, 24: 44--68, 2018.

    \bibitem{zhou2008}
     X.~Zhou, W.~Ren and W.~E.
    \newblock {Adaptive minimum action method for the study of rare events}.
    \newblock {\em The Journal of Chemical Physics}, 128: 104--111, 2008.

    \bibitem{weinan2004}
     W.~E, W.~Ren and E.~Vanden-Eijnden.
    \newblock {Minimum action method for the study of rare events}.
    \newblock {\em Communications on Pure and Applied Mathematics}, 57: 637--656, 2004.

    \bibitem{wan2018}
     X.~Wan, H.~Yu and J.~Zhai.
    \newblock {Convergence analysis of a finite element approximation of minimum action methods}.
    \newblock {\em SIAM Journal on Numerical Analysis}, 56: 1597--1620, 2018.  
    
\end{thebibliography}
\end{document}